\documentstyle[12pt]{article}
\begin{document}
\newtheorem{lem}{Lemma}
\newtheorem{rem}{Remark}
\newtheorem{question}{Question}
\newtheorem{prop}{Proposition}
\newtheorem{cor}{Corollary}
\newtheorem{thm}{Theorem}
\newtheorem{conj}{Conjecture}

\title
{Partial results on extending
the Hopf Lemma }
\author{YanYan Li\thanks{Partially
 supported by
      NSF grant DMS-0701545.}
\\
Department of Mathematics\\
Rutgers University\\
110 Frelinghuysen Road\\
Piscataway, NJ 08854\\
\\
Louis Nirenberg
 \\
Courant Institute\\
251 Mercer Street \\
New York, NY 10012\\
}

\date{ }
\maketitle

Dedicated with affection to Umberto Mosco on his 70'th birthday

\input { amssym.def}

\bigskip

\bigskip

\setcounter {section} {0}

\section{Introduction}

In \cite{LN1}, Theorem \ref{thm2.1},
the authors proved, in one dimension, a generalization of the Hopf
Lemma, and the question arose if it could be extended to higher
dimensions.  In this paper we present
two conjectures as possible extensions, and give
a very partial
answer.  We write this paper to call attention to the problem.

The one dimensional result of  \cite{LN1} was
\begin{thm}
Let $u\ge v$ be \underline{positive}
$C^3$, $C^2$
functions 
respectively on $(0, b)$
which are also in $C^1([0, b])$.  Assume
\begin{equation}
u(0)=\dot u(0)=0
\label{1.1}
\end{equation}
and
$$
\mbox{either}\ \dot u>0
\ \mbox{on}\ (0, b)\
\mbox{or}\ \dot v>0\ \mbox{on}\ (0,b).
$$
\underline{Main condition:}
\begin{equation}
\mbox{whenever}\
u(t)=v(s)\
\mbox{for}\ 0<t\le s<b,
\ \mbox{there}\ 
\ddot u(t) \le v''(s),
\label{1.2}
\end{equation}
(here $\cdot=\frac{d}{dt}$,
$' = \frac d{ds})$.

Then
\begin{equation}
u\equiv v\ \ \mbox{on}\
[0, b].
\label{1.3}
\end{equation}
\label{thm1.1}
\end{thm}

The proof given in \cite{LN1} is somewhat
roundabout.
In
the Appendix we present a more direct one, but it is still a bit 
tricky.
In \cite{LN1}, it was assumed that $u$ is of class $C^2$ on
$(0, b)$, but its proof there 
actually required that $u$ be of class $C^3$.

\bigskip

\bigskip

Turn now to higher dimensions.  
Let  $u\ge v$
be $C^\infty$ functions of $(t, y)$,
$y\in \Bbb R^n$, in 
$$
\Omega=\{(t, y)\ |\
0< t< 1, |y|<1\},
$$
and $C^\infty$ in the closure of $\Omega$.
Assume that
\begin{equation}
u>0, \ v>0,\ u_t>0\quad \mbox{in}\ \Omega
\label{1.4}
\end{equation}
and
\begin{equation}
u(0, y)=0\quad\mbox{for}\ |y|<1.
\label{1.5}
\end{equation}

We impose a main condition:
\begin{equation}
\mbox{whenever}\ u(t,y)=v(s,y)\
\mbox{for}\ 0<t\le s<1, |y|<1,  
\mbox{there}\ \Delta u(t,y)\le \Delta v(s,y).
\label{1.6}
\end{equation}
Under some additional conditions  we wish to conclude that
\begin{equation}
u\equiv v.
\label{1.7}
\end{equation}

Here are two conjectures, in decreasing strength,
which would extend Theorem \ref{thm1.1}.  In each,
we consider $u$ and $v$ as above.

\begin{conj} \label{conjecture1}
Assume, in addition, that
\begin{equation}
u_t(0, 0)=0.
\label{1.8}
\end{equation}
Then (\ref{1.3}) holds:
$$
u\equiv v.
$$
\end{conj}

\begin{conj} \label{conjecture2}
In addition to (\ref{1.8}) assume that
\begin{equation}
u(t, 0)\ \mbox{and}\ v(t, 0)\
\mbox{vanish at}\ t=0\
\mbox{of finite order}.
\label{1.9}
\end{equation}
Then 
$$
u\equiv v.
$$
\end{conj}

We have not succeeded in proving them.
What we present here is a  partial answer
to Conjecture \ref{conjecture2}:  Here let
$k, l$ be the orders of the first $t-$derivative of $u$, $v$ respectively
at the origin which are not zero.
Clearly $k\le l$.

\begin{thm}\label{thm1.2}
In addition to  the conditions of Conjecture \ref{conjecture2}, we assume
the annoying condition 
\begin{equation}
\nabla_yu_{tt}(0, 0)=0.
\label{extra}
\end{equation}
Then $
u\equiv v$ provided $k=2$ or $3$.
\end{thm}

For $k<3$ the proof is simple, but not that
for $k=3$.

We will always use Taylor series expansions for $u, v$, in $t$,
\begin{equation}\label{1.12}
u=a_1(y)t+a_2(y)\frac {t^2}{2!}
+a_3(y)\frac {t^3}{3!}+\cdots,\qquad 
v=b_1(y)t+b_2(y)\frac {t^2}{2!}
+b_3(y)\frac {t^3}{3!}+\cdots
\end{equation}
The conditions on $u$ and $v$ are as follows
\begin{equation}\label{1.13}
0\le u(t)-v(t)
=(a_1-b_1)t
+(a_2-b_2) \frac {t^2}{2!}
+(a_3-b_3) \frac {t^3}{3!}+\cdots
\end{equation}
where
$$
u(t, y)=v(s, y),\quad t\le s,
$$
i.e.
\begin{equation}\label{1.14}
a_1(y)t +a_2(y)\frac {t^2}{2!}
+a_3(y)\frac {t^3}{3!}+\cdots
=b_1(y)s+b_2(y)\frac {s^2}{2!}
+b_3(y)\frac {s^3}{3!}+\cdots,
\end{equation}
there
\begin{equation}\label{1.15}
0\ge \Delta u-\Delta v
=(a_2-b_2)+
t(\Delta a_1+a_3)-s(\Delta b_1+b_3)
+ \frac {t^2}2
(\Delta a_2+a_4)
- \frac {s^2}2
(\Delta b_2+b_4)+\cdots
\end{equation}

We first present the proof of the more difficult case
$k=3$.  It takes up sections 2-5.
In section 6 we treat the case $k=2$.

\section{}

\noindent{\bf Steps of the proof.}\
We are assuming $k=3$.  
The proof consists of two steps:

\noindent{\bf  Step A.}\
This consists in proving
\begin{thm}\label{thm2.1}
Under the conditions of Theorem \ref{thm1.2},
where $k=3$, we have
\begin{equation}\label{2.1}
l=3, \mbox{and}\  b_3(0)=a_3(0).
\end{equation}
\end{thm}

\noindent{\bf Step B.}\
In this step we consider our condition
\begin{equation}\label{2.2}
u(t, y)=v(s, y)\ \mbox{for}\ 
0\le t\le s.
\end{equation}
Since $u_t>0$ for $t>0$, we may solve this for
$t=t(s,y)$.  Assuming that $u$ is not identically
equal to $v$, for 
\begin{equation}\label{2.3}
\tau(s,y) =s-t(s, y)
\end{equation}
we derive, from (\ref{1.6}), an
elliptic  
differential inequality  for $\tau(s, y)$.  Using a comparison function
we prove that 
\begin{equation}\label{2.4}
\tau(s, 0)\ge \epsilon s\ \ 
\mbox{for some}\ 0<\epsilon\
\mbox{small}.
\end{equation}

On the other hand, for $y=0$, we have,
by (\ref{2.1})  and (\ref{1.12}),
$$
u(t, 0)=v(s, 0)
$$
i.e. after dividing by $a_3(0)$,
$$
t^3+\mbox{higher order terms}=s^3+
\mbox{higher order terms}.
$$
Hence 
$$
t(s, 0)=s+
\mbox{higher order terms}.
$$
But this contradicts (\ref{2.4}), and the proof of Theorem \ref{thm1.2}
is then complete.

For $k=3$, we will first present the
proof of Step B; it seems more interesting to us.

\section{} 
{\bf
  Proof of (\ref{2.4}) in case $k=3$.
}

Here we assume that (\ref{2.1}) holds, i.e.
$$
b_3(0)=a_3(0)>0
$$
and 
first derive the elliptic inequality for
$
\tau(s, y)=s-t(s,y),
$ where $t(s,y)$ is the solution of 
\begin{equation}\label{3.1}
u(t(s,y), y)=v(s, y).
\end{equation}

Differentiating this we find, setting $v_i=
\partial_{y_i}v$,
$$
v_s=u_t t_s,\qquad v_{ss}=u_t t_{ss}+u_{tt}t_s^2,
$$
$$
v_i =u_t t_i +u_i,\qquad
v_{ii}=u_t t_{ii}+ 2u_{ti}t_i+u_{tt}t_i^2+u_{ii}.
$$
Hence 
$$
0\le \Delta v(s,y)-\Delta u(t,y)
=u_t \Delta t+2u_{ti}t_i+u_{tt}
(|\nabla t|^2-1).
$$
In terms of $\tau=s-t$, this becomes, after dividing by $u_t$,
\begin{equation}
F(\tau):=\Delta \tau -\frac{u_{tt}}
{u_t} (|\nabla \tau|^2-2\tau_s)+2 \frac {u_{ti}}{u_t}\tau_i\le 0.
\label{3.2}
\end{equation}

This is the differential inequality for $\tau$. 

We will consider this in the region
\begin{equation}
D=\{(s, y)\ |\ s>K|y|^2\},
\quad K\ \mbox{large}, \
\mbox{near the origin},
\label{3.3}
\end{equation}
and use a comparison function:
\begin{equation}
\label{3.4}
h=s+s^{1+\delta}-C|y|^2,\quad
\delta=\frac 14, C=K+1.
\end{equation}
Near the origin we have
\begin{equation}\label{3.5}
h(s,y)\le 0\ \ \mbox{where}\
s=K|y|^2.
\end{equation}

We assume now that
$v$ is not identically equal to $u$ near the origin and
argue by contradiction.

Observe first that if $v(\bar s, \bar y)=
u(\bar s, \bar y)$ for some $\bar y$ and
some $\bar s>0$ then
$\tau (\bar s, \bar y)=0$.
But near $(\bar s, \bar y)$,
$\tau\ge 0$ satisfies the inequality (\ref{3.2}),
which is elliptic there.
By the strong maximum principle, it would
follow that $\tau\equiv 0$ there.
Then, again by the strong maximum principle
$\tau\equiv 0$ everywhere,
i.e. $v\equiv u$ near the
origin, for $t\ge 0$.  Contradiction.

Thus we may assume that
$\tau>0$ for $s>0$

The basic result of this section is
\begin{lem}\label{lem3.1}  For $0<\epsilon$ small,
$\tau\ge \epsilon h$ in $D$ near the origin.
\end{lem}

Once the lemma is proved, it follows that
$\tau(s, 0)\ge \epsilon s$ for $0<s$ small,
i.e., (\ref{2.4}) holds, and Step B would be complete.

\noindent {\bf Proof of Lemma \ref{lem3.1}}.\
Choose positive $\epsilon\le 1/10$, so small
that on $D\cap \{s=c\}$, $c$ to be fixed --- where
 $\tau$ is positive, and hence bounded away from
zero ---
\begin{equation}
\tau\ge \epsilon h,
\label{3.6}
\end{equation}
$\epsilon$ depends on $c$.

In view of (\ref{3.5}) it follows then
that
$$
\tau-\epsilon h\ge 0, 
\mbox{on the boundary of 
}\ G=D\cap \{s<c\}.
$$

We now use the maximum principle, suitably to show that

\begin{equation}\label{3.7}
\tau\ge \epsilon h\quad\mbox{in}\ G.
\end{equation}
--- completing the proof of Lemma \ref{lem3.1}.
We argue by contradiction.

Suppose $\tau-\epsilon h$ has a negative minimum at
some point $(\bar s, \bar y)$ in $G$.  There,
of course, 
$$
\tau<\epsilon(s+s^{1+\delta})<2\epsilon s,
$$
and so
\begin{equation}\label{3.8}
t=s-\tau\ge (1-2\epsilon)s\ge \frac 45 s.
\end{equation}

At $(\bar s, \bar y)$, $\nabla \tau =\epsilon \nabla h$ and
$$
\Delta \tau\ge \epsilon \Delta h.
$$
Therefore, there, $\epsilon h$ satisfies the inequality
$$
\Delta (\epsilon h) - 
\frac{u_{tt}}
{u_t}  (\epsilon^2 |\nabla h|^2 -2
\epsilon h_s) +2\epsilon \frac{  u_{ti}  }{ u_t} h_i
\le 0
$$
i.e.  after dividing by $\epsilon$,
\begin{equation}\label{3.9}
 F[\epsilon h]=
\Delta h - 
\frac{u_{tt}}
{u_t} \left\{ \epsilon[( 1+(1+\delta)s^\delta)^2 +4C^2|y|^2]
-2- 2(1+\delta)s^\delta\right\}
-4C \frac {u_{ti} y_i}  {u_t}
\le 0.
\end{equation}

For small $\epsilon$ and $c$ (which may depend on $K$), 
\begin{equation}\label{3.11}
\mbox{
the expression \{\ \  \} in
(\ref{3.9}) is negative}.
\end{equation}

We will choose $K$ to ensure that
\begin{equation}\label{3.12}
u_{tt}(t(\bar s, \bar y), \bar y)\ge 0.
\end{equation}
We have
\begin{equation}\label{3.13}
u_{tt}=a_2+a_3 t+\cdots
\end{equation}
Since $a_3(0)>0$, near the origin,
\begin{equation}\label{3.14}
a_3(t, y)\ge \frac {a_3(0)}2.
\end{equation}
Recall that $u_t>0$, i.e.
\begin{equation}\label{3.15}
0<a_1+ta_2+ \frac{t^2}2 a_3+\cdots
\end{equation}
Thus $a_1\ge 0$ and $a_1=O(|y|^2)$.  
By (\ref{extra}), and it is only here that 
(\ref{extra}) is used, 
\begin{equation}\label{3.16}
|a_2|\le A|y|^2
\end{equation}
for some $A>0$.

Now, still at $(\bar s, \bar y)$, and for $t=t(\bar s, \bar y)$,
we have
\begin{eqnarray*}
u_{tt}&=&
a_2+a_3t+\cdots\ge \frac{  a_3(0) }2 t-A|y|^2+O(t^2)\\
&\ge & \frac {a_3(0)}4 t-A|y|^2\qquad \quad (\mbox{for}\ c\ \mbox{small})\\
&\ge & \frac {a_3(0)}5s
-A|y|^2
\end{eqnarray*}
by (\ref{3.8}).  We require
$$
K\ge \frac {5A}{ a_3(0)}.
$$
Then (\ref{3.12}) holds:
$$
u_{tt}\ge 0.
$$
(we may suppose $K>1$.)

Consequently, from (\ref{3.9}) we find
\begin{equation}\label{3.17}
\Delta h-\frac {4C}{u_t}
u_{ti}y_i\le 0\qquad
\mbox{at}\
(\bar s, \bar y).
\end{equation}

Next, by a well known elementary inequality
which uses the fact that the second order derivatives in $y$ of $u_t$ are bounded in absolute value we have, for some constant $B$,
$$
|u_{ti}|\le B\sqrt{ u_t}\qquad \forall\ i.
$$
So
\begin{equation}\label{3.18}
M:= \frac{4C}{ u_t}
|u_{ti} y_i|\le
\frac{4CB|y|}{  \sqrt{u_t}  }.
\end{equation}

Now, recall, $t=t(\bar s, \bar y)$,
$$
u_t=a_1+a_2 t+\frac {a_3t^2}2+\cdots
\ge t(a_2+ \frac {a_3t}2+\cdots)
\ge t(-A|y|^2 + \frac {a_3(0)}4 t)
$$
by (\ref{3.16}), for $t$ small.  So
$$
u_t\ge t(-\frac AK s+  \frac {a_3(0)}4 t)
\ge \frac 45 s( -\frac {A}Ks +  \frac {a_3(0)}5 s)
$$
by (\ref{3.8}).  Hence
\begin{equation}\label{3.19}
u_t \ge \frac{a_3(0)}{10}s^2
\end{equation}
provided
\begin{equation}\label{3.20}
\frac AK \le \frac{a_3(0)}{100}.
\end{equation}

Inserting (\ref{3.19}) in (\ref{3.18}) we find
\begin{equation}\label{3.21}
M=\big|\frac{4C}{ u_t} \sum u_{ti}y_i|\le
\frac{L|y|}s
\end{equation}
where
$$
L= \frac{4\sqrt{10}CB}{ \sqrt{a_3(0) } }.
$$

Thus, by (\ref{3.3}),
$$
M\le \frac L{ \sqrt{Ks} }.
$$

We now insert this in (\ref{3.17}) and, computing
$\Delta h$, we find
$$
\delta(1+\delta) s^{\delta-1}-2n C\le 
\frac { 4\sqrt{10}  }{  \sqrt{a_3(0)}  }
\frac{  K+1}{ \sqrt{K}  }\frac B{  \sqrt{s} }.
$$
But for $\delta =1/4$, and $c$ restricted still further
if necessary, we see that this is impossible.

\vskip 5pt
\hfill $\Box$
\vskip 5pt

\begin{rem}\label{rem1}  Our use of the maximum principle is somewhat 
unusual.   Normally, one
would prove that $F[\epsilon h]$, in
(\ref{3.9}) is positive in
$G$; in fact we do not know how
to prove that.  But, as we see,
 it suffices only to show that it is positive
at $(t(\bar s, \bar y), \bar y)$.
\end{rem}

\section{Step A}

\noindent{\bf 4.1.}\  We turn now to Step A.  Let
\begin{equation}\label{4.1}
\hat a_i(y)\ \mbox{be the lowest order terms of}\
a_i(y)
\end{equation}
in its Taylor expansion; $\hat a_i$ is a homogeneous
polynomial.  We know that
\begin{equation}\label{4.2}
\deg \hat a_1, \deg \hat b_1, \deg (\hat a_2-\hat b_2)\ge 2,
\end{equation}
since, by (\ref{1.15}), $\hat a_2-\hat b_2$ is non-positive.

Our aim is to prove, in this and the next
section, that
if $k=3$ then
\begin{equation}\label{4.3}
l=3\ \ \mbox{and}\
b_3(0)=a_3(0).
\end{equation}

We will constantly use (\ref{1.13})-(\ref{1.15}).

\noindent{\bf Proof that if }\ $l=3$ {\bf then}\
$b_3(0)=a_3(0)$.

Since $u\ge v>0$ in $\Omega$, necessarily
$$
a_3(0)\ge b_3(0)>0.
$$

In (\ref{1.14}) set $y=0$ and solve for
$t=t(s)$.  Clearly
$$
t=  \left(  \frac{b_3(0)}{ a_3(0) }\right)^{\frac 13} s
+O(s^2).
$$
Inserting this value for $t(s)$ in (\ref{1.15})
we find, by looking
at the coefficients,
$$
0\ge \left( \frac{b_3(0)}{a_3(0)}\right)^{ \frac 13} (\Delta \hat a_1(0)+
a_3(0))
-(\Delta \hat b_1(0)+b_3(0)).
$$
i.e.
\begin{equation}\label{4.10}
(b_3)^{ \frac 13}  \Delta \hat a_1 -
(a_3)^{ \frac 13} \Delta \hat b_1
+(b_3)^{ \frac 13} a_3
-(a_3)^{ \frac 13} b_3\le 0,\qquad\mbox{at}\ y=0.
\end{equation}

Since $a_3\ge b_3>0$ at $y=0$, we infer that
\begin{equation}\label{4.11}
(b_3)^{ \frac 13}  \Delta \hat a_1 -(a_3)^{ \frac 13} \Delta \hat  b_1
\le 0,\qquad\mbox{at}\ y=0.
\end{equation}

Now $\hat a_1\ge \hat b_1\ge 0$. This implies
$\Delta \hat a_1(0)\ge \Delta \hat  b_1(0)\ge 0$.
  If both $=0$ then
(\ref{4.10}) implies $a_3(0)=b_3(0)$.

Then, since $\Delta \hat a_1(0)>0$, it follows that
\begin{equation}
\Delta \hat b_1(0)> 0.
\end{equation}
In particular, $\deg \hat b_1=\deg \hat a_1=2$.

Next, at a point $y$ where $\hat b_1(y)> 0$, take
$$
s=K\hat a_1(y), \ K\ \mbox{large}.
$$
Then from (\ref{1.14}) we solve for $t=t(s)$ and
find, looking at terms of various degrees in $y$,
$$
t= K\hat b_1(y) +\circ(|y|^2).
$$

Insert this in (\ref{1.15}); we obtain, looking at terms of
second degree in $y$, and using the fact that $K$ is arbitrarily large,
\begin{equation}\label{4.13}
0\ge \hat b_1(y)(\Delta \hat a_1(0)+a_3(0))-
\hat a_1(y) (\Delta \hat b_1(0) +b_3(0)).
\end{equation}

Since the right hand side is a homogeneous quadratic, its
Laplacian is $\le 0$, i.e.
$$
0\ge \Delta \hat b_1(\Delta \hat a_1 +a_3(0))-\Delta \hat a_1
(\Delta \hat b_1+b_3(0)),
$$
so
$$
a_3(0)\Delta \hat b_1-b_3(0)\Delta \hat a_1\le 0.
$$
Using (\ref{4.11}) it follows, then, that
$$
a_3^{ \frac 23}
b_3^{\frac 13} \Delta \hat a_1\le b_3\Delta \hat a_1
$$
which implies (\ref{4.3}):
$$
b_3(0)=a_3(0).
$$

\bigskip

From now on we assume $l>3$ and prove that this is
impossible.

\noindent{\bf 4.2.}\ {\bf The
 case 
 $l>3$}.

{\bf (i)}\ {\bf Claim 1}\ In this case 
\begin{equation}\label{4.14}
b_1=O(|y|^4).
\end{equation}

\noindent{\bf Proof.}\ Suppose not, then $\hat b_1$ has degree $2$ since
by 
 the positivity of $v$, 
$\hat b_1\ge 0$. $\hat a_1$ also has degree $2$ since $a_1\ge b_1$. 
The proof above of (\ref{4.13}) still works, and yields
\begin{equation}\label{4.15}
0\ge \hat b_1(\Delta \hat a_1+a_3(0))-\hat a_1 \Delta \hat b_1.
\end{equation}

Taking trace we find
$$
0\ge \Delta \hat b_1 a_3(0)
$$
i.e. $\hat b_1=0$ --- recall that $\hat b_1\ge 0$.  Contradiction.
The claim is proved.

Next, set $y=0$ and solve for $t(s)$ in (\ref{1.14}).  We find
$$
t=\left(  \frac 6{l!} \frac {b_l(0)}  { a_3(0)  }\right)^{1/3} s^{l/3}
+\circ(s^{l/3}).
$$
Inserting this in (\ref{1.15}) we find, at $y=0$, since
$\Delta \hat b_1=0$,
$$
0\ge 
\left(\frac 6{l!}\right)^{\frac 13}
\left(\frac{b_l}{a_3}\right)^{1/3}s^{l/3} (\Delta a_1+a_3)-s^2
(\Delta b_2+b_4)+\circ(s^{l/3}+s^2).
$$
Consequently
$$
l\ge 6.
$$

We shall make use of the following
\begin{lem} \label{lemma2}
Let $v\ge 0$ be given by (\ref{1.12}) and assume that
 $l$ is the order of the first $t-$derivative of
$v$ which is $>0$ at the origin.  Let $m$ be the
first value of $i$ (if it exists) such that
$$
\deg \hat b_{i}=1.
$$
Suppose that for some $j$, $1\le j\le (l+4)/3$,

$$
\deg \hat b_i\ge 3\ \ \mbox{for}\ i<j.
$$
Then
\begin{equation}
 m\ge \frac {l+j}2.
\label{52}
\end{equation}
\end{lem}

\noindent{\bf Proof.}\ Clearly $j\le m<l$.   At some $y$,
$\hat b_m(y)<0$.  Then, at that $y$, if we set
$$
s=|y|^a, \qquad 0<a\ \ \mbox{to be chosen},
$$
we have, since $v\ge 0$,
\begin{equation}
0\le
\sum_{i< j}\frac 1{i!}b_{i}(y)s^{i}
+ \sum_{j\le i\le m-1}\frac 1{i!} b_{i}(y)s^{i}
+ \sum_{m\le i\le l-1} \frac 1{i!} b_{i}(y)s^{i}
+O(s^l).
\label{48prime}
\end{equation}

In case $j=1$ we find
\begin{equation}
0\le -\frac 1{2m!}\hat b_m s^m
=O(|y|^2s)+O(s^l).
\label{48primeprime}
\end{equation}

Suppose that (\ref{52}) does not hold, i.e.
$$
m<\frac {l+1}2.
$$
Then there exists $a>0$ such that $\deg$ LHS of
(\ref{48primeprime})
$<$  $\deg$ of each term on RHS of (\ref{48primeprime}).
One easily verifies this
using the fact that
$$
\frac 1{l-m}<\frac 1{m-1}.
$$
But then (\ref{48primeprime})
is impossible.

In case $j>1$ we find from (\ref{48prime})
and the fact that $\hat b_1=O(|y|^4)$, that 
\begin{equation}\label{53}
0
\le -\frac { \hat b_m(y) |y|^{am} }{ 2m!}
\le O(|y|^{4+a})+ 
O(|y|^{3+2a})+
O(|y|^{2+ja})+O( |y|^{la}).
\end{equation}

Suppose that (\ref{52}) does not hold, i.e.
\begin{equation}\label{54}
m<\frac{l+j}2.
\end{equation}

\noindent{\bf Claim:}\ There exists $a>0$ such that 
the degree of LHS of (\ref{53}) 
$<$ the degree of each term on RHS of (\ref{53}).

If so, (\ref{54}) is impossible.

\noindent{\bf  Proof of Claim.}\
The claim asserts the existence of $a>0$
such that
\begin{equation}
\label{55}
\left\{
\begin{array}{rl}
1+ma<4+a, \ \ \mbox{i.e.}\ a<\frac 3{m-1},\\
1+ma<3+2a, \ \ \mbox{i.e.}\ a<\frac 2{m-2}\ \mbox{if}\
m>2,\\
1+ma<2+ja, \ \ \mbox{i.e.}\ a< \frac 1{m-j}\ \mbox{if}\
m>j,\\
1+ma<la,  \ \ \mbox{i.e.}\ a>\frac 1{l-m}.
\end{array}
\right.
\end{equation}

If  $m=2$,  the second and third inequalities automatically hold, so does
the third if $m=j$.  Otherwise it says that
$$
a<\frac 1{m-j}.
$$

One easily verifies using (\ref{54})
that
$$
\frac 1{l-m}
< 
\left\{
\begin{array}{rl}
 \frac 3{m-1}, &
\mbox{if}\ m=j=2,
\\
\min\{ \frac 3{m-1},  \frac 2{m-2}\}&
\mbox{if}\ m=j
\ge 3,\\
\min\{ \frac 3{m-1}, \frac 2{m-2}, \frac 1{m-j}\}&
\mbox{if}\ m>j.
\end{array}
\right.
$$
It follows that the required $a$ exists.  Hence, Lemma 
\ref{lemma2} is proved. 

\section{}
We come now to  
a crucial step.

\begin{prop}  If $l\ge 3i$, $l>3$, $i\ge 1$, then
$$
\deg \hat b_i\ge 3.
$$
\label{proposition1}
\end{prop}

Using the proposition we may now
give the

\noindent{\bf Completion of the proof of Theorem \ref{thm2.1}.}\
At $y=0$, if we solve (\ref{1.14})
for $t$ we find as before,
$$
t=As^{l/3}+\circ(s^{l/3}),
$$
where
$$
A=(\frac 6{l!} \frac {b_l}{a_3})^{1/3}.
$$
Inserting this in (\ref{1.15}) and using Proposition
\ref{proposition1} we see that
$$
0\ge As^{l/3}(\Delta a_1+a_3)+O(s^{  [l/3]+1  }).
$$
But this is impossible, and Theorem \ref{thm2.1} is proved.

\noindent{\bf Proof of Proposition \ref{proposition1}.}\ 
By Lemma \ref{lemma2},
$$
\deg \hat b_i>1\quad \mbox{for}\ i<\frac l2+1.
$$
Suppose the proposition is false.  Then there is a first $j\le l/3$ such that
$$
\deg \hat b_j=2.
$$
We will show that this is impossible.

By (\ref{4.14}), $j\ge 2$.

\noindent{\bf Claim.}\ $\hat b_j\ge 0$.

If not, at some $y$, $\hat b_j(y)<0$.  Then, setting
$$
s=|y|^a,
$$
we have, using Lemma \ref{lemma2}, and (\ref{4.14}), 
\begin{equation}
0<-\frac{\hat b_j |y|^{ja}}{2j!}
=O(|y|^{4+a})+O(|y|^{2a+3})
+O(|y|^{ 1+a(l+j)/2})
+O(|y|^{al}).
\label{56}
\end{equation}
Setting $a>1/j$ but very close to $1/j$,
 we see that
the degree in $y$ of LHS of (\ref{56})
$<$ the degree of each term on RHS of (\ref{56}), i.e.
(here we use $j\le l/3$)
\begin{equation}
2+ja<
\min\{4+a, 2a+3,  1+a(l+j)/2, al\}.
\label{sss}
\end{equation}

But then (\ref{56}) is impossible.  The claim is proved.

We now distinguish two cases.

\noindent{\bf Case 1.}\ $\deg \hat a_1=2$.  We have $\hat b_j\ge 0$.

Fix $y$ so that $\hat b_j(y)>0$; since $\hat a_1$
cannot vanish on an open set we may also ensure that $\hat a_1(y)>0$.

As before, set $s=|y|^a$, with $a>1/j$ but very close to $1/j$,
so that (\ref{sss}) holds.  Then, as before,
in the expression for $v$ the term
\begin{equation}
J=
\frac 1{j!} \hat b_j(y)s^j=
\frac 1{j!} \hat b_j(y)|y|^{aj}
\label{ttt}
\end{equation}
has degree smaller than that of any other term.

Consequently we may solve (\ref{1.14}) first, and find 
$$
t=\frac { \hat b_j(y) }{ j! \hat a_1(y)} 
|y|^{aj}+\circ(|y|^{aj}).
$$
Inserting these values for $s$ and $t$ in (\ref{1.15})
we find
$$
0\ge \frac {|y|^{aj}} {j!}
\frac{  \hat b_j}{ \hat a_1} (\Delta\hat  a_1 +a_3(0))
-\frac {|y|^{aj}}{ j!} \Delta \hat b_j+\circ(|y|^{aj}),
$$
i.e.
$$
0\ge \hat b_j (\Delta \hat a_1 +a_3(0))
-\hat a_1 \Delta \hat b_j.
$$
As before, taking trace, we conclude that
$\hat b_j=0$.  Contradiction.

\noindent{\bf Case 2.}\
$\deg \hat a_1>2$.  Then $\deg \hat a_1\ge 4$.

Still take 
 $s=|y|^a$, with $a>1/j$ but very close to $1/j$,
so that (\ref{sss}) holds. 
We still have that in the expression for $v$, the term
$
J$ in (\ref{ttt})
 has degree smaller than that of
every other term.  To  solve (\ref{1.14}) for $t$, we note that
the leading terms of $u(t, y)$ are now
$$
u(t, y)=a_1(y)t +\frac 12 a_2(y)t^2+
\frac 16 a_3(y)t^3+\cdots
=O(|y|^4t)+O(|y|^2t^2)+a_3(0)t^3+\cdots,
$$
where we have used $\deg \hat a_2\ge 2$ which follows from
Lemma \ref{lemma2}.
Thus
$$
t=\left( \frac 6{ a_3(0) }J\right)^{\frac 13} +
\circ(|y|^{ \frac {2+aj}{3} }).
$$
Inserting 
 these values
for $s$ and $t$ in (\ref{1.15}) we find
$$
0\ge t a_3(0)-\frac {s^j}{j!} \Delta \hat b_j
+\circ(|y|^{ \frac {2+aj}{3} })+\circ(|y|^{aj}).
$$
It follows, since $(2+aj)/3<aj$,
that
$0\ge a_3(0)$,
a
 contradiction.

The proof of Proposition \ref{proposition1} in
case $\deg \hat a_1>2$ is complete.
Theorem \ref{thm2.1} is proved.
\vskip 5pt
\hfill $\Box$
\vskip 5pt

\section{ Proof of Theorem \ref{thm1.2} in case $k=2$}

The proof has again Step A and Step B.
i.e. we first prove that
\begin{equation}\label{6.1}
l=2\ \mbox{and}\ b_2(0)=a_2(0),
\end{equation}
and then if $u$ is not identically equal to $v$,
using the differential inequality (\ref{3.2})
for $\tau$, and the same comparison function
$h$ of (\ref{3.4}) we derive a contradiction.

The proof of (\ref{6.1}) is trivial: from (\ref{1.13}), 
$$
a_2(0)-b_2(0)\ge 0
$$ while
from (\ref{1.15}), at $t=0$, the opposite inequality holds.

Turn now to the equation for $\tau$.  We follow the argument of section 3.  We have
to prove that $\tau-\epsilon h$ cannot
have a negative minimum in $G$.  To do this
we have to check, as before that
$F[\epsilon h]$ in (\ref{3.9}) is
positive at a possible 
minimum point
$(\bar s, \bar y)$, i.e.
\begin{equation}\label{6.2}
\delta(1+\delta) \bar s^{-\delta -1}-2nC -
\frac {u_{tt}}{ u_t}
\{\quad\}- \frac {4C u_{ti}\bar y_i}{ u_t}>0.
\end{equation}
The term $\{\quad\}<0$, and $u_{tt}=
a_2+O(t)>0$, since
$a_2(0)>0$.  In addition,
$$
M=\frac {4C}{ u_t}
|u_{ti}\bar y_i|\le
\frac { 4C \sqrt{ \sum |u_{ti}|^2 } |\bar y|}{ u_t}.
$$
Now 
$$
u_t=a_1+ a_2t+\cdots
\ge \frac 12 a_2(0)t>\frac 25 a_2(0) s
$$
by (\ref{3.8}).  Thus, since
$s>K|y|^2$,
$$
M\le \frac {10C |\nabla^2 u|}{ a_2(0)\sqrt{K}
\sqrt{s} }.
$$
We conclude that (recall $C=K+1$),
$$
 F[\epsilon h]
\ge \delta(1+\delta)s^{\delta -1}-2nC -\mbox{constant} \cdot
\frac{ \sqrt{K}}
{\sqrt{s}}>0
$$
since $\delta=1/4$.  (\ref{4.2})
is proved, and the proof of Theorem \ref{thm1.2} for $k=2$ is complete.

\vskip 5pt
\hfill $\Box$
\vskip 5pt

\section{Appendix.  A simple proof of Theorem \ref{thm1.1}}

We treat only the case:
\begin{equation}
\dot u>0\quad\mbox{on}\ (0, b).
\label{A.1}
\end{equation}

We have to prove that
\begin{equation}
u\equiv v.
\label{A.2}
\end{equation}

The proof proceeds in two steps:

\noindent{\bf Step A.}\ (\ref{A.2}) holds in case
\begin{equation}
v'(s)\ge 0.
\label{A.3}
\end{equation}

\noindent{\bf Step B.}\ Necessarily,
$$
v'(s)\ge 0.
$$

{\bf Step A.}\
Proof of (\ref{A.2}) if $v'\ge 0$.

We have
$$
u(t)=v(s),
$$
since $u'>0$, for $t>0$, we may
solve for $t=t(s)$.
Here $\cdot=
\frac d{dt}$,
$'=\frac d{ds}$.
Then 
$$
v'=\dot u t'.
$$
Compute
\begin{eqnarray}
(v'^2 -\dot u^2)'&=&
2v'v'' -2 \dot u \ddot u t'
=2v' (v''-\ddot u)\label{A.4}\\
&\ge& 0
\nonumber
\end{eqnarray}
by our main condition (\ref{1.2}).  But at the origin,
$$
v'^2 -\dot u^2=0,
$$
so 
$$
v'^2 -\dot u^2=
\dot u^2(t'^2 -1)\ge 0.
$$
Hence 
$$
t'^2\ge 1.
$$
Since $t'\ge 0$ somewhere
for $s$ arbitrarily small, it follows
that $t'\ge 1$, i.e. $t\ge s$.  But
then $t\equiv s$ and so $u\equiv v$.

\vskip 5pt
\hfill $\Box$
\vskip 5pt

{\bf Step B.}\ Proof that $v'\ge 0$.

(i)\ We use part of an argument of \cite{LN1}:
$$
\ddot u(t)\ \mbox{is a function of}\ t
$$
but since $\dot u>0$ it may be written as a function of
$u$, i.e.
\begin{equation}
\ddot u=f(u),
\label{A.5}
\end{equation}
with, however, $f$ an unknown function.
$f$ is continuous on an interval
$[0, m]$ for some $m>0$, and of class $C^1$ on
$(0, m]$, since $u$ is of class $C^3$ for $t>0$.

The main condition (\ref{1.2}): 
$$
\ddot u(t)\le v''(s)\quad \mbox{whenever}\
u(t)=v(s),\ t\le s,
$$
is equivalent to the inequality
\begin{equation}
\label{A.6}
v''\ge f(v).
\end{equation}

We have $u\ge v$ and both vanish, with their first 
derivatives at the origin.  But we cannot apply the Hopf Lemma
to $(u-v)$ because $f$ is not known to be 
Lipschitz near the origin.

\begin{lem} \label{lemA.1}
If $v(s)=u(s)$ for some
$s>0$, then
$$
v\equiv u.
$$
\end{lem}

\noindent {\bf Proof.}\  We use a differential inequality which 
holds for $\tau=s-t(s)$.  Namely, we have
$$
v'=\dot u t',
$$
$$
v''= \dot u t'' +\ddot u t'^2
=-\dot u \tau'' +\ddot u(1-\tau')^2.
$$
So
$$
0\le v'' -\ddot u= -\dot u \tau''  +\ddot u(\tau'^2-2\tau').
$$

Now if $u(s)=v(s)$ for some $s>0$, then,
there, $\tau =0$.  But $\tau\le 0$.  By the strong
maximum principle
it would follow that $\tau\equiv 0$, i.e. $v\equiv u$.

\vskip 5pt
\hfill $\Box$
\vskip 5pt

To prove that $v'\ge 0$ we argue by contradiction.  Suppose $v'<0$ somewhere.

(ii)\ We cannot have 
$v' \ge 0$ on an interval $(0, c)$, for
if this holds, by Step A, we would have
$$
v\equiv u\quad\mbox{on}\ (0, c).
$$
By Lemma \ref{lemA.1}, we would have
$$
v\equiv u\quad\mbox{everywhere}.
$$

So, arbitrarily near the origin 
there are points where $v'<0$.
But then there must be an interval
$(a, c)$, $0<a<c<b$
on which 
$$
v'<0\ \mbox{and}\ v'(a)=0.
$$
On this interval, by (\ref{A.4}),
$$
(v'^2 -\dot u^2)'\le 0.
$$
Hence 
$$
v'(s)^2 -\dot u(t(s))^2 \le -\dot u^2(t(a))\qquad\mbox{on}\
(a,c)
$$
and, consequently,
$$
\dot u(t(a))\le \dot u(t(s))\qquad\mbox{for }\ a<s<c.
$$
It follows that
$$
\ddot u(t(a))\ge  0.
$$
By our main condition, then
$$
v''(a)\ge \ddot u(t(a))\ge 0.
$$

Now we cannot have 
$v''(a)>0$ since 
$0=\dot v (a)>\dot v(s)
$ for $a<s<c$.  Thus 
\begin{equation}\label{A.7}
v''(a)=0,\quad 
\mbox{and so}\ \ddot u(t(a))=0.
\end{equation}

(iii)\ We now make use of 
(\ref{A.5}) and (\ref{A.6}).
By (\ref{A.5}), 
$$
0=f(u(t(a)))=f(v(a)).
$$
Hence, by (\ref{A.6}), on $(a, c)$,
$$
v''(s)\ge f(v(s))=f(v(s))-f(v(a))
=f'(\xi)(v(s)-v(a))
$$
for some $\xi$ in $(v(s), v(a))$.

But $v(s)-v(a)$ has its maximum at
$a$.  We may apply the classical Hopf Lemma to
infer that
$$
v'(a)<0.
$$
This contradicts the fact that $v'(a)=0$.

\vskip 5pt
\hfill $\Box$

\end{document}